# Transforming Post-Secondary Education in Mathematics

*Driving constructive change in education in the mathematical sciences at two-year colleges, four-year colleges, and universities*

Tara Holm, Cornell University

## Abstract

In this manuscript, I introduce and describe the work of mathematicians and mathematics educators in the group Transforming Post-Secondary Education in Mathematics (TPSE Math or TPSE, for short). TPSE aims to coördinate and drive constructive change in education in the mathematical sciences at two-year colleges, four-year colleges, and universities across the nation. It seeks to build on the successes of the entire mathematical sciences community.

This manuscript reviews the events that led to the founding of TPSE Math and articulates its vision and mission. In its first phase with national events, TPSE found broad consensus with the mathematical sciences community on the challenges facing the community. Learning from educational transformations experiences in other scientific fields, and with the support of the Mathematical Advisory Group of 34 mathematical sciences department chairs and leaders, TPSE moves into a second phase focused on action. This is a snapshot in time, and TPSE's ongoing activities will continue to be documented and disseminated. The piece concludes with a reflection of the impact that my involvement in this work has had on my career.

## 1 Introduction

The education landscape has changed dramatically in the last half century. Higher education has become essential to economic mobility. At the same time, colleges, universities, and students are under severe financial pressure. And new pedagogies and technologies make it possible to reach students in many more ways. These and other forces will change higher education (Bok 2013, Chapter 10).

Mathematics departments play a central role in undergraduate education: few departments teach a larger percentage of the student body. Mathematicians must respond to the challenges facing higher education. If we opt out, we risk losing the substantial role that mathematics departments currently play, and we endanger the health of the US mathematical sciences research enterprise.

A volunteer group of mathematicians and mathematics educators called **Transforming Post-Secondary Education in Mathematics** (TPSE Math or TPSE, pronounced "tipsy", for short) is working to support the mathematics community in this endeavor. Eric Friedlander, President of the American Mathematical Society (AMS) 2011–2012, invited me to join the TPSE leadership team when I was Chair of the AMS Committee on Education (CoE). In what follows, I describe the history leading to the current state of post-secondary education in mathematics. I then chronicle the formation of TPSE and the foundation of work that it builds on. I report on TPSE's current partnerships and plans. I conclude with my personal history and perspective on TPSE.

## 2 Landscape preceding TPSE

In the past twenty years, there have been many calls to improve mathematics instruction. With particular attention to research universities, the AMS Task Force on Excellence exhorted, "To ensure their institution's commitment to excellence in mathematics research, doctoral departments must pursue excellence in their instructional programs" (Ewing 1999, p. 3). Departments must maintain a relevant and broad curriculum. In addition to what they must teach, departments must also address questions of how to teach it (Ewing 1999). In 2003, Halpern and Hakel (2003, p. 38) noted, "it would be difficult to design an educational model that is more at odds with the findings of current research about human cognition than the one being used today at most colleges and universities."

The AMS Task Force on First-Year Mathematics (Lewis & Tucker 2009, p. 755) made three key suggestions towards the pursuit of excellent instruction:

1. Harness the power of technology to improve teaching and learning;
2. Leadership matters – success in this area depends upon the value assigned to it by a department's leadership; and
3. Invest in teaching graduate students to be good teachers.

Preparing graduate students to teach, a particular role for research universities, is intertwined with any discussion of undergraduate education: graduate students represent the future of the professoriate. This highlights the particular role research universities play in our profession.

In their landmark Proceedings of the National Academy of Sciences study, Freeman et al. (2014) leave no doubt that active learning techniques improve student performance. Doctoral departments must adjust graduate teaching preparation accordingly, and offer renewed support to all mathematics faculty members. Guidance for graduate students and current faculty must be broad, as no one technique will work in every classroom.

There are several aspects of the current post-secondary education landscape that particularly inform TPSE Math choices. Every five years, the Conference Board of the Mathematical Sciences (CBMS) undertakes a statistical study of undergraduate programs in mathematics. One trend that their data make clear is that more and more students are enrolling in courses at two-year colleges. In 2010, over 45% of enrollments in mathematics, statistics and computer science courses taught in mathematics departments were taught in two-year colleges (Blair et al. 2013, Table S.1). As a consequence of this, TPSE made the commitment early on to include two-year colleges and questions of credit transfer as an integral part of its work. Students in the 21$^{st}$ century are different from their mid-20$^{th}$ century counterparts. Babcock and Marks (2011) have analyzed a number of datasets to report that from 1961 to 2004, the amount of time a typical undergraduate spends on academic work has dropped by nearly a third, from around 40 hours per week to 27. This is no doubt influenced by the dramatic increase in the cost of post-secondary education. The financial pressures that both students and universities face affect the types of programs that are feasible. Simultaneously, the increasing costs ratchet up the need to articulate the value mathematics adds to post-secondary education. Finally, there are now very good techniques available to assess the effectiveness of educational initiatives. TPSE Math is committed to evaluating the success of the projects it pursues and the innovations it endorses. It aims to cultivate collaborations between mathematicians, researchers in mathematics education and evaluators to strengthen assessment procedures.

There is renewed federal interest in higher education in general, and undergraduate science, technology, engineering and mathematics (STEM) education in particular. President Obama identified post-secondary education as key to a stronger economy and 21st century success of the nation. He asked the President's Council of Advisors on Science and Technology (PCAST) to prepare a report on producing one million more STEM graduates over the next decade. In that report, PCAST points to a US Department of Commerce report that projects a 17% increase in the need for STEM-trained graduates over this time period (PCAST 2012). The mathematics community was taken aback when PCAST suggested that "faculty from mathematics-intensive disciplines other than mathematics" should develop and teach courses in college-level mathematics, and that there should be a "new pathway for producing K–12 mathematics teachers from … programs in mathematics-intensive fields other than mathematics" (PCAST 2012, p. 30). When writing its report, PCAST did not broadly consult the mathematics community. They were surprised not to find a journal focused on undergraduate mathematics education among AMS journals. They were not aware of the many successful innovations in post-secondary mathematics instruction (P. LePage, personal communication, October 22, 2015). This is an important lesson for the mathematics community: we must redouble efforts to promulgate our successes beyond just our community.

In a more positive light, the National Research Council has described how mathematics has become essential to modern science, and recommends that undergraduate education in the mathematical sciences reflect this new stature (Everhart et al. 2013). The good news is that all of these recommendations have spurred action within the mathematics community.

Calls for innovation and transformation of mathematics instruction are not new. Still the current spotlight on mathematics comes at a time when there is broad awareness of the challenges we face and an increased focus on educational outcomes. Moreover, the mathematical sciences community has been and is increasingly involved in developing solutions. This is a once-in-a-generation opportunity to capitalize on the power of collective action and support the transformation of post-secondary education in the mathematical sciences.

## 3 Formation of TPSE Math

In February 2013, Carnegie Corporation of New York assembled a group of higher education leaders in the mathematical sciences to take stock and envision how to transform the field from a service discipline to an essential partner in post-secondary education. As a result of this meeting, Phillip Griffiths founded the group TPSE Math. TPSE has listened to and will continue to work with the mathematical sciences community to determine how best to achieve systemic change. It is now poised to forge alliances with and among state and federal agencies, the policy community, university administrators, higher education associations, and professional organizations to secure the financial and structural support necessary to achieve these goals.

In May 2016, TPSE Math incorporated as an Educational Program affiliated with the University System of Maryland Foundation. The Board of Governors includes the following.

1. **Phillip Griffiths** (Board Chair) is Professor Emeritus of Mathematics and former Director of the Institute for Advanced Study. He was Provost of Duke University.
2. **Eric Friedlander** is the Dean's Professor of Mathematics at the University of Southern California and is a Past President of the American Mathematical Society (AMS).
3. **S. James Gates, Jr.** is the John S. Toll Professor of Physics at the University of Maryland and a member of the President's Council of Advisors on Science and Technology (PCAST).
4. **Mark Green** is Professor Emeritus of Mathematics at the University of California, Los Angeles, and former Director of the Institute for Pure and Applied Mathematics (IPAM).
5. **Tara Holm** is Professor of Mathematics at Cornell University and former Chair of the AMS Committee on Education.
6. **Karen Saxe** is the DeWitt Wallace Professor of Mathematics at Macalester College, a past Vice President of the Mathematics Association of America (MAA), and a leader of MAA's Common Vision Project.
7. **Uri Treisman** is Professor of Mathematics and of Public Affairs at the University of Texas, Austin. He founded and directs the Charles A. Dana Center.

In 2015, **William (Brit) Kirwan**, Chancellor Emeritus of the University of Maryland, joined TPSE as a Senior Advisor. He chairs the Conference Board of the Mathematical Sciences (CBMS) and was appointed Executive Director of TPSE Math in May 2016.

TPSE Math envisions a future where postsecondary mathematics education will enable any student, regardless of his or her chosen program of study, to develop the mathematical knowledge and skills necessary for productive engagement in society and in the workplace. TPSE's mission statement articulates its vision (TPSE Math 2015):

> *TPSE Math will facilitate an inclusive movement to strengthen postsecondary education in mathematics by working closely with – and mobilizing when necessary – faculty leaders, university administrations, membership associations, and relevant disciplinary societies in the pursuit of mathematically rich and relevant education for all students, whatever their chosen field of study. TPSE Math will identify innovative practices where they exist, advocate for innovation where they do not, and work with and through partners to implement and scale effective practices.*

TPSE takes a national-level approach to transformation, seeking to leverage resources from non-profit foundations and federal agencies to increase the capacity of the profession to achieve change. The mathematical sciences community must proceed with coherence but not uniformity, ever heeding local needs of individual institutions to ensure appropriate changes will take root.

From 2013 to 2015, TPSE surveyed what is happening in the mathematics community through one national and four regional meetings. It engaged strategy consultants from Parthenon-EY to help evaluate the state of post-secondary education in the mathematical sciences. Through our meetings and research, TPSE has found a high level of consensus among faculty and administrators about the need for renewal of the post-secondary curriculum. Professional development opportunities can support interested faculty to develop and enhance their pedagogical practice. TPSE evaluated the capacity of the professional societies to support community-wide change. The following have been identified as areas where the TPSE leadership can have the biggest impact, by working in concert with existing programs and societies to leverage existing capacity towards a common goal:

1. **Curriculum pathways** (lower and upper division, allowing students to reach the mathematics relevant to their field of study);
2. **Graduate co-curricular training**; and
3. **Leadership development**.

The TPSE Math leadership team has been successful at engaging large swathes of the mathematics community. Through national and regional meetings, it has learned from and fostered conversations among departmental leaders from a broad range of institutions, college and university administrators, and representatives of federal and non-profit funding agencies. TPSE has brought together leaders of those professional societies that include post-secondary mathematics education as part of their primary mission:

1. American Mathematical Association of Two-Year Colleges (AMATYC)
2. American Mathematical Society (AMS)
3. American Statistical Association (ASA)
4. Association for Women in Mathematics (AWM)
5. Conference Board of the Mathematical Sciences (CBMS)
6. Mathematical Association of America (MAA)
7. Society for Industrial and Applied Mathematics (SIAM)

Starting in 2016, TPSE turned to the action phase of its work.

## 4 Foundation on which TPSE builds

### 4.1 History of change in US education in mathematics

The last period of dramatic change in high school and college mathematics curricula began in the 1950s (Tucker 2015). In 1957, the Soviets launched Sputnik I, the first satellite to go into orbit. Four months later, the US launched its first successful satellite, Explorer I. The Cold War Space Race had begun. In that era of unprecedented public support for science education, calculus became the ultimate goal of high school mathematics. Supported in part by the Ford Foundation, AP calculus came into being. Since then, the exam has shifted to being a test of calculus knowledge rather than more general problem solving (Bressoud et al. 2012). The variety of mathematics relevant to the world has expanded remarkably in the 60 years since then. We must open new pathways to offer students the mathematics they need. This is a particular challenge in mathematics, where theories do not become false or go out of fashion. We need to find new trails through the mountain range of the mathematical sciences.

### 4.2 Work by the Professional Societies

One role that the AMS CoE plays is to coöperate with the other professional mathematics societies on matters concerning education. At the October 2015 AMS CoE meeting, three other societies gave presentations about their current projects. David Bressoud (Macalester College) provided an analysis of the MAA's studies of first-semester calculus instruction across in the US (Bressoud 2015). Donna Lalonde (ASA) reported on the varied outreach efforts that ASA undertakes to engage students at all levels in the statistical sciences (LaLonde 2015). Rachel Levy (Harvey Mudd College) is working with SIAM to create infrastructure to support internships in business, industry and government for mathematics students (Levy 2015). SIAM has also reported on the value of mathematical modeling in the K-12 and undergraduate curriculum (Turner 2014), reflecting the increasing use and applications of mathematics across disciplines.

The MAA has been a leader of the mathematical community on the topic of undergraduate mathematics and a tireless supporter of educational initiatives. Roughly every 10 years since 1953, the Committee on the Undergraduate Program in Mathematics (CUPM) has produced the benchmark curriculum guide for mathematics departments (MAA 2015). While the other professional societies are consulted during the preparation of each guide, the MAA provides the lion's share of support for the process. There is a new effort underway to create a companion instructional practices guide that will provide would-be innovators and experts alike with a basic framework for developing new, evidence-based classroom practices. CUPM also has a subcommittee focused on Curriculum Renewal Across the First Two Years (CRAFTY) whose 2011 report (MAA 2011) is particularly relevant to TPSE's work. One of the broadest impact the MAA has had may be through the New Experiences in Teaching Project (Project NExT). Since 1994, this professional development program has helped new faculty members become more effective teachers. After the yearlong program ends, NExT Fellows have access to the NExT listserv, providing continuous support and discussion about all matters loosely related to teaching.

More recently, the MAA coördinated the *Common Vision* project, bringing together representatives from the five professional societies concerned with post-secondary mathematics education – AMATYC, AMS, ASA, MAA, SIAM. Three members of the Common Vision leadership team are also affiliated with TPSE Math (the principal investigator Karen Saxe, Uri Treisman and myself). This project included a workshop and culminated in a report (Saxe & Braddy 2015) detailing the commonalities among the five societies' curricular recommendations and recording the shared opinions of the pressing need to transform post-secondary education in the mathematical sciences (Holm & Saxe 2016).

## 4.3 Innovations in Doctoral Departments

Each October, the AMS CoE meets in Washington DC and hosts a forum for discussion of issues in mathematics education. During my tenure as Chair of the CoE, the focus of the Committee's work shifted from K–12 education to post-secondary mathematics education, as this is the aspect of mathematics education most closely related to the daily work of AMS members. While we have come to expect a high level of commitment to undergraduate education in mathematics departments at primarily undergraduate teaching institutions, the CoE has seen an impressive array of programs at research universities. A few examples from 2014 and 2015 meetings are described below, and more details are available at the AMS CoE website (AMS 2016).

Both the University of Michigan and its Mathematics Department are deeply committed to their educational mission. The 1999 AMS Task Force on Excellence found "a culture in the mathematics department that encourages and rewards innovation, one that is well rounded, that strikes a balance between teaching and research, and that supports the work of students and colleagues at all levels" (Ewing 1999, p. 84). In his 2014 presentation, Stephen DeBacker of the Mathematics Department reported that innovations in active learning in first-year calculus have taken root and have now expanded to labs for second-year courses. Preparing graduate students and postdoctoral scholars to teach in the Michigan program ensures consistency and success. He cited as essential the full support of the university administration in this endeavor, and concluded, "A successful undergraduate program requires the efforts of nearly everyone in the department" (DeBacker 2014).

The Mathematics Department at the University of Illinois has developed several successful programs that reach a wide range of students studying mathematics (Ando 2014). Mathematicians have collaborated with engineering faculty to develop workshops where students must apply their calculus knowledge to solve problems inspired by real-world applications in engineering. The University uses Treisman's collaborative learning model in its Merit Program to support students from traditionally underrepresented populations. Merit scholars are more successful than their peers in Illinois calculus. Finally, the Illinois Geometry Lab provides research experiences for 40 undergraduates each semester.

Through his San Francisco State University–Colombia Combinatorics initiative, Federico Ardila (2015) has forged connections between US students and their Colombian peers, each serving as role models for the other. Starting with the two principles that mathematical ability is uniformly distributed and that every student can have a meaningful mathematical experience, Ardila took an intentional approach to building a bridge between the two communities, using technology to offer courses and facilitate groups working together between the two countries. Of the 200 students who have participated in the program, 45 are currently working on their doctorates, including 26 women and 20 students from underrepresented groups.

The CoE also heard about the important role that the EDGE program and the National Alliance for Doctoral Studies in Mathematics have played in mentoring the next generation of mathematical scholars, ensuring a diverse talent pool (Wilson 2015; Math Alliance 2013).

## 5 Building Systemic Change

In spite of the successes detailed here and the many more described on the AMS CoE web site, we still face significant challenges. Very few efforts are scaled or transferred; many rely on a charismatic individual for their continuation. Large research universities play a significant role in setting the standards for post-secondary education, but the mathematics research community has not been engaged in a coördinated way with undergraduate mathematics education. There is little more than informal coöperation among the mathematics professional societies, not through lack of will but simply because each society is focused on its own mission. CBMS, the umbrella organization of all mathematical sciences societies, is understaffed and has too broad a mission to spearhead an initiative on undergraduate education. Despite the good will of its members, it is not in a position to generate change at the national scale.

By contrast, TPSE Math has a very precise focus, and strong connections to college and university administrations, higher education associations, foundations, federal agencies, state government offices, and mathematics departments and societies. We aim to reinforce and augment the relationships that the professional societies have already built in the policy sphere. For these reasons, TPSE is especially well positioned to forge and sustain effective partnerships that will propel transformation in post-secondary education in mathematics.

Generating systemic change is a notoriously complex challenge. Fortunately, there are models that have been successful in academia and can be adapted for the mathematical sciences community. In the Life Sciences and in Physics, curricula and pedagogy are now better adjusted to foster students' conceptual understanding of the science. In both cases, real progress occurred only after the communities came together to articulate a coherent vision. Disconnected innovations are insufficient to transform the entire field.

In mathematics there are successful programs such as the ones described above that we can build upon to facilitate successful propagation of change. Still, mathematics is different from the other sciences. In the physical sciences, where research is dependent on expensive equipment and experiments, the community decides with funding agencies on the top research priorities for the field. As a consequence, the scientists are more accustomed to work together as a community. An important TPSE goal is to enhance the existing structures within the mathematical sciences community to ensure this necessary community-wide progress.

## 5.1 Life Sciences

The biological sciences changed dramatically in the second half of the 20[th] century. The ultimate goal of understanding life remains constant, but the types of questions life scientists can ask and the tools available to answer them have developed rapidly. Beginning in 2008, the American Association for the Advancement of Science (AAAS) started a series of conversations with faculty, administrators, students, biological sciences professional societies, and funding agencies on the future of undergraduate biology education. First-year courses were completely overhauled, with careful attention to the desired outcomes for biology students and for general education students. Using the *Vision & Change in Undergraduate Biology Education* recommendations as a guide (Vision & Change 2010), the Partnership for Undergraduate Life Sciences Education (PULSE) has developed tools that promote department-wide implementation of new curricula and evidence-based pedagogies (PULSE 2014. A key feature of PULSE is its use of social connections, both in workshops and through a faculty ambassador network, to propagate change. This has been shown to lead to greater adoption of innovations than evidence presented via literature (Henderson & Dancy 2011). The landscape in the biological sciences is far more diverse than in the mathematical sciences. For example, there are dozens of professional societies of scientists in the life sciences. Nevertheless the essential aspects of PULSE could transfer to the mathematical sciences community.

## 5.2 Physical Sciences

Physicists and physics educators have received much attention for the development of a number of pedagogical techniques and assessment tools. The Force Concept Inventory (FCI) was the first test of its kind (Hestenes et al. 1992). It measures the change in students' conceptual understanding from the beginning to the end of a first physics course. Jerome Epstein has developed a Calculus Concept Inventory, which has not been adopted to nearly the same extent as the FCI (Epstein 2007). Eric Mazur and his colleagues have demonstrated the power of peer instruction to increase students' conceptual learning (Crouch & Mazur 2001). Carl Wieman was awarded the National Science Foundation (NSF) Director's Award for Distinguished Teaching Scholars in 2001, the same year that he won the Nobel prize. Wieman strongly advocates for activities that keep students actively engaged in exploring physics during class time. There are significant differences between first-year physics classrooms and first-year mathematics ones. In physics, students may be more homogeneous and there is not nearly as much variation in what a student's first post-secondary physics course could be. Still, we should analyze how and why it is that new instructional strategies are so much more quickly implemented by physical scientists, while being mindful that the reality of innovation may not be as widespread as self-reports lead us to believe (Dancy & Henderson 2010). It is also important to remember that active engagement in doing mathematics is a long-held value of many mathematicians. R. L. Moore and his followers have used a version of inquiry-based learning since the first half of the 20[th] century (Wilder 1976). This technique is used in many advanced mathematics courses, and more recently is being adapted for appropriate use in first-year courses.

## 5.3 First steps in the Mathematical Sciences

The mathematics community has achieved some transformation at the state level: mathematics departments at all post-secondary public institutions in individual states are beginning to work together. In 2013, the Ohio Mathematics Steering Committee was charged with "develop[ing] expectations and processes that result in each campus offering pathways in mathematics that yield

1. Increased success for students in the study of mathematics;
2. A higher percentage of students completing degree programs; and
3. Effective transferability of credits for students moving from one institution to another." (Ohio Mathematics Initiative 2014, p. 2)

This is indicative of how degree completion now dictates state policy in higher education. Ohio's college and university mathematics department chairs identified key challenges in achieving these goals and are now making progress towards addressing them. More and more states are following suit with their own mathematics task forces, many with the help of my TPSE Math colleague Uri Treisman and the Dana Center. Going forward, transfer of credits from one higher education institution to another will become an increasing challenge, especially given financial pressures that are causing more students to divide their education among multiple institutions.

## 6 Work by TPSE Math

From the beginning, TPSE has been intent on engaging the entire mathematics community. It has organized events at the Joint Mathematics Meetings (JMM), as well as a national and a series of regional TPSE meetings. The goal has been and is to work with the mathematics community to identify the most urgent issues, see how they are being addressed, and determine which early experiments and models can be scaled up and used by others. TPSE has also sought opinions from the "demand" side of the equation, including employers and disciplinary partners.

## 6.1 Preparation

TPSE's first public event was a standing-room-only panel discussion at the Baltimore JMM in 2014 (TPSE Math 2014a). Moderated by Philip Griffiths, panel members included Michelle Cahill (Carnegie Corporation of New York), Jo Handelsman (Yale University and now Associate Director for Science at the White House Office for Science and Technology Policy), Brit Kirwan and Joan Leitzel (Ohio Mathematics Initiative). Jerry McNerney, the only US Congressman with a PhD in Mathematics, also joined the conversation, providing valuable insights from his perspective on Capitol Hill. The discussions clarified how current governmental policy affects higher education in general and mathematics in particular. All agreed that we have a moment of opportunity for innovation and transformation of post-secondary mathematics education.

## 6.2 Information gathering

In June 2014, TPSE organized a national meeting at the University of Texas, Austin, to bring together leaders from academia, business and government to discuss challenges and explore scalable solutions (TPSE Math 2014b). It was at this meeting that the following questions were articulated:

1. How can the **undergraduate curriculum** be reshaped to raise the level of numeracy among citizens and better align current teaching with the expanded role of mathematics?
2. How can **new pathways** be opened to enable non-majors and developmental students to reach the level of numeracy needed for careers that demand analytical thinking and 21st century quantitative skills?
3. How will **new technologies and teaching** trends affect pedagogy and the economic model of mathematics departments?
4. How can a broader, more relevant undergraduate experience better prepare students for the **workplace of the future**, including interdisciplinary opportunities?
5. How can **graduate students** be equipped to teach more broadly about the uses of mathematics while maintaining depth in their own research?

Next came four regional meetings. Each meeting attracted participants from a variety of colleges, universities, and funding agencies. The meetings were open to all, but TPSE began by inviting participants identified through local contacts and the professional societies to include a diverse group of stakeholders. We structured the regional meetings around panel discussions on a wide variety of topics. These are summarized as follows.

1. In November 2014, at the University of Maryland Baltimore County (UMBC), TPSE arranged three panels: (1) a description of TPSE Math and other initiatives with similar goals, (2) disparities in participation by various populations, and (3) issues facing non-R1 institutions. There were 50 participants.
2. In February 2015, at the University of California, Los Angeles (UCLA), TPSE organized three panels: (1) the role of mathematics in career preparation, (2) the role of (the field of) mathematics education in post-secondary education in the mathematical sciences, and (3) system-wide efforts to improve post-secondary education in the mathematical sciences. There were 49 participants.
3. In September 2015, at the University of Chicago, TPSE organized a two-day meeting with six panels: (1) the role of college and university administrations, (2) the role of mathematics departments, (3) preparation of graduate students as future faculty, (4) secondary school teacher training, (5) enhanced opportunities for highly motivated undergraduates, and (6) the role, relevance and reform of calculus. There were 63 participants.
4. In December 2015, at Duke University, TPSE organized a two-day meeting with five panels: (1) multiple math pathways, (2) math and other disciplines, (3) math courses for non-STEM undergraduates, (4) adaptive learning, and (5) statistics and big data. At this meeting, TPSE announced some of its next steps in the action phase (Section 6.4 below). There were 66 participants.

Further details and videos of some panels are available on the Meetings page of the TPSE Math website (TPSE Math 2015).

While each regional meeting had a slightly different focus, there emerged from all meetings a broad consensus among the mathematics community that all departments are facing pressure on these issues, some of the important work to address them has begun, and that we as a community must find a coherent way forward that allows for local variation. Moreover, transformation cannot rely on continuous additional resources. The declining state support of public universities is unlikely to reverse. Departments may be able to get one-time allocations to support a phase transition to more innovative curricula and pedagogies. To succeed in the long term, sustainability needs to be built in from the beginning.

## 6.3 Building relationships

In addition to organizing meetings, the TPSE leadership has been actively engaged in building relationships with the professional societies, department leaders and mathematicians more broadly. It is committed to intentionally building a community of mathematics leaders who collectively drive the transformation process. At the San Antonio 2015 JMM, TPSE organized a discussion with the leadership of the professional societies and associations. It was an unusual gathering, possibly the first time that the leadership of AMATYC, AMS, the Association of Public and Land-grant Universities (APLU), CBMS, MAA, the National Council of Teachers of Mathematics (NCTM), and SIAM, as well as directors from NSF, all met at the same time and place. In March 2015, I had the opportunity to engage with mathematics department chairs from research universities at MSRI's Sponsors' Day event. Mark Green, Uri Treisman and I gave an update on TPSE Math's activities to CBMS at their May 2015 meeting (Green et al.

2015). At all of these meetings, there was a strong consensus that now is the time to promote transformation.

## 6.4 Action Phase

Grounded in the consensus it has found, TPSE Math now aims to advance several projects to achieve transformation of post-secondary education in the mathematical sciences. TPSE has identified the areas where it is most likely to have a significant impact. These include partnering with associations and societies to achieve shared goals and seeking to provide a scaffolding to support the mathematics community efforts.

At a meeting in Washington DC in March 2015, TPSE initiated the Mathematics Advisory Group (MAG) (TPSE Math 2016). Funded by the Sloan Foundation, we have identified 34 department chairs and leaders have been identified who will gather to begin to develop an action plan to carry out, scale up and evaluate the effectiveness of major reforms. TPSE sought leaders who represent the gamut of higher education institutions, and who bring a diversity of views and experiences. This core group will be a key action and communication partner, advising TPSE on "grass roots" issues at the departmental level and helping identify successful and valuable models. The MAG will help TPSE convene a larger meeting of 100 to 200 department chairs and leaders to share information and establish partnerships of departments committed to reform in post-secondary education in the mathematical sciences. Members of the MAG may become TPSE Math Ambassadors, willing to advise departments about transformation. This is the key first step in intentionally developing the leadership and increasing the capacity of the mathematics community. Through the MAG and Ambassador network, TPSE also plans to reinvigorate the preparation of graduate students as educators and mentors.

TPSE Math has also embarked on several partnerships to enhance curriculum pathways. It will serve as an advisory partner to ITHAKA S+R (Ithaka 2004), a research group that studies the use of technology to improve teaching and learning, as well as the economic impact of such technologies. TPSE will also serve as an advisory partner to APLU, the American Association of State Colleges and Universities (AASCU), and the Dana Center in developing multiple pathways in lower division mathematics courses to improve completion rates and quality of instruction. TPSE also plans to promote renewal and creation of upper division curricula in response to the growing demands from other disciplines.

With the recent appointment of Kirwan as Executive Director, TPSE Math is preparing for the next stage of action, seeking substantial financial support from a number of non-profit donors and federal agencies. This action phase also signals the beginning of sustained efforts to increase the diversity of the mathematics community, and to build expertise in data analytics to better assess needs and evaluate outcomes.

## 7 A personal journey to TPSE Math

As a faculty member in a research-intensive mathematics department, my principal job is mathematics research. I am also a practitioner in mathematics education, and I strive to use research findings in mathematics education to inform my classroom practice. Through my service in the AMS and now with TPSE Math, I work on transforming undergraduate mathematics education at the national level.

## 7.1 Path to involvement

My personal interest in undergraduate mathematics education dates back to my days as a high school student taking my first college mathematics course. I am fortunate to have had outstanding teachers, and I am indebted to my mentors at each stage of my education.

As a graduate student at MIT, given only a couple days of TA training, my initial teaching assignments were recitations where I was forbidden to lecture. Instead I was supposed to engage students in problem solving. I was one of three math graduate students to take an education course at MIT, discussing practical and theoretical aspects of general undergraduate education. That course certainly highlighted for me some of the ways in which mathematics is similar to and yet very different from other disciplines. I had the good fortune to spend a couple months teaching mathematics to computer science students in a non-traditional program at the short-lived ArsDigita University (ArsDigita 2002). This was an opportunity to teach incredibly motivated students using problem sets guiding them through parts of calculus, statistics and discrete mathematics. Together with Shai Simonson (Stonehill College), we engaged the students in a new-to-them research project that involved mathematics and computer experimentation to understand a card trick (Holm & Simonson 2003).

After graduate school, I spent three years at UC Berkeley, funded in part by an NSF postdoctoral fellowship. I taught three courses during my time there, and experienced first-hand the challenge of teaching upper division courses in large lecture format. After a year at the University of Connecticut, I started a tenure-track position at Cornell University. I was happy to find a position in a strong research department that also has a deep commitment to excellent undergraduate education. I was the first Cornell faculty member to participate in MAA's Project NExT. The sessions and materials from NExT workshops and the ongoing support through the electronic network have proved a tremendous resource in my teaching. One particular session, Joe Gallian's (University of Minnesota, Duluth) advice to "just say Yes" to opportunities to give back to the profession, has influenced my approach to service at my university and through the professional societies.

Soon after earning tenure, I was asked to run for election to the governing boards of the AMS and the AWM; I was elected to each in 2011 and 2012 respectively. AMS Council members also serve on one of the five policy committees: Education; Meetings and Conferences; the Profession; Publications; or Science Policy. I was assigned to the Committee on Education (CoE), and chaired it from 2012 to 2016. During that time period, Eric Friedlander was President of AMS. One of his primary objectives as President was for the Society to increase its participation in the improvement of post-secondary education (E. M. Friedlander, personal communication, 29 October 2011). We worked together to shift the AMS CoE focus from K-12 to post-secondary education in the mathematical sciences. My fortuitous appointment to the CoE has opened opportunities to have an impact on undergraduate education in the mathematical sciences in ways that I could never have predicted. In particular, I was invited to join the leadership teams for the Common Vision project and TPSE Math as a direct result of my work with the AMS CoE.

The most pleasant surprise in my role as Chair of the AMS Committee on Education (CoE) was the opportunity to laud the work that mathematicians and mathematics educators are doing. In my AMS and more general mathematical travels, every mathematics department I have visited or heard about has some interesting project afoot. Faculty members want their students to engage deeply in thinking about mathematics. The professional societies seek to support their members' teaching mission. The mathematical sciences community must acknowledge and publicize these existing successes. Identifying the most promising innovations and determining the best way to adapt them and scale them for use at different institutions is no small task; this is at the heart of TPSE's mission.

A skeptic might point out that AMS CoE and TPSE Math meetings are assemblies of the willing, and ask whether the broader mathematics community is on board. Anecdotally, I have found general interest and support from mathematicians in the mathematics research community. For example, at the Spring 2016 Texas Geometry and Topology Conference, I was asked to give both a research talk and a second talk about my work with the AMS CoE and TPSE Math. Both talks were well attended, and the latter generated thoughtful discussions among all participants, from graduate students to the most senior topologists in the room (Holm 2016). Moving forward, when TPSE identifies departments to serve as lodestars, it will be important to select those that have high levels of faculty commitment to transformation. Their early successes will serve as models to promote change at all institutions.

## 7.2 Reflections

I conclude with a personal perspective. As indicated earlier, I am not a mathematics education researcher; but I believe strongly that the mathematics community must maintain the bridges between researchers in mathematics education and practitioners of mathematics education, particularly at the post-secondary level. Moving forward we need to improve our communication and collaboration. Working with the AMS CoE and with TPSE Math, I have had the tremendous opportunity to engage with the senior leadership in mathematics and policy leaders in academia and government.

I count myself lucky to be a member of a supportive research department where faculty members are encouraged to contribute to all aspects of the profession. I have no illusions: my work with the AMS and TPSE Math did not get me tenure or promotion to full professor. It was considered a favorable part of my dossier, but my research is the *sine qua non*. These service opportunities did arise at a good time in my career. I had young children at home. Particularly while my second child was an infant and I was on parental leave, I appreciated the opportunity to be engaged with the mathematics community in this way. TPSE was in an early phase when most of the work consisted of phone meetings and email correspondence. This all fit into to the spare time I might find at odd times of the day. TPSE now involves more travel, but it is work that I continue to be able to fit in with the rest of my research and teaching. I welcome it as a chance to think at the community-wide and national level about the future of our profession.

Cornell University does provide strong support for its faculty. The Office of Faculty Development and Diversity offers a number of professional development opportunities and mentoring programs. For example, they have encouraged faculty to raise their voices beyond the walls of academe by offering a Public Voices fellowships through the Op-Ed project. The Center for Teaching Excellence offers workshops, lunches and logistical support for faculty who want to innovate in their classrooms (Cornell 2012). They supported the Mathematics Department in bringing the Discovering the Art of Mathematics leaders (Fleron et al. 2008) to Cornell to sponsor a workshop for Cornell and Ithaca College faculty members introducing their teaching materials for introductory general education courses and their inquiry-based approach. Implementing innovative teaching practices has also been an attractive cause for university fund raising. By the end of 2015, donations had funded nearly $1million in grants to faculty members to support curricular renewal and innovation through the newly formed Engaged Cornell program (Cornell 2016).

Through my work with TPSE Math, I have come to understand better the political and financial forces that are reshaping the way the public and the mathematics community perceive the role of mathematics in today's society and for the future. All mathematics departments are under pressure: the federal and state governments are curtailing their contributions to universities; university administrations are slashing resources; and everyone is demanding more from higher education institutions. The alignment of these forces sets the stage for great change in the teaching and learning of mathematics. It is my hope that all levels of the mathematics community – from department colleagues and administrators to the leadership of professional societies – will come together to ensure that our students are prepared for a future we cannot yet imagine.

**Acknowledgements.** Tara Holm thanks the anonymous referees and the editors for detailed reports that greatly improved the structure of this Chapter. TPSE Math is deeply grateful for the continuing support of the Arthur P. Sloan Foundation and Carnegie Corporation of New York.